\documentclass{svjour3}

\newcommand{\norm}[1]{\Vert#1\Vert}

\newcommand{\Real}{\mathbb R}

\renewcommand{\u}{\mathbf{u}}
\renewcommand{\v}{\mathbf{v}}

\renewcommand{\r}{\mathbf{r}}

\newcommand{\w}{\mathbf{w}}
\newcommand{\Wb}{{\boldsymbol{W}}}

\newcommand{\phib}{{\boldsymbol{\phi}}}
\newcommand{\Phib}{{\boldsymbol{\Phi}}}
\newcommand{\alphab}{{\boldsymbol{\alpha}}}
\newcommand{\gammab}{{\boldsymbol{\gamma}}}

\newcommand{\Rc}{\mathcal{R}}
\newcommand{\Tc}{\mathcal{T}}
\newcommand{\Sc}{\mathcal{S}}

\newcommand{\Uc}{\mathcal{U}}

\newcommand{\Fc}{\mathcal{F}}

\newcommand{\Pbb}{\mathbb{P}}
\newcommand{\Rbb}{\mathbb{R}}
\newcommand{\Wbb}{\mathbb{W}}
\DeclareSymbolFontAlphabet{\mathcal}{symbols}

  \usepackage{pgfplots}
  \pgfplotsset{compat=newest}
  \usetikzlibrary{plotmarks}
  \usetikzlibrary{arrows.meta}
  \usepgfplotslibrary{patchplots}
  \usepackage{grffile}
  
\usepackage{graphicx}
\usepackage{color}
\usepackage{epsfig}
\usepackage{bm}
\usepackage{xspace}
\usepackage{amsmath}
\usepackage{amssymb}
\usepackage{amsbsy}
\usepackage{eucal}
\usepackage{mathrsfs}
\usepackage{fancyvrb}
\usepackage{color, colortbl}
\usepackage{subfigure}
\usepackage{multirow}
\usepackage{algorithm}
\usepackage{algorithmic}
\usepackage{makecell}

\definecolor{LightCyan}{rgb}{0.88,1,1}

\newcommand{\jb}{\bm {j}}

\newcommand{\Mb}{\bm {M}}

\newcommand{\Sb}{\bm {S}}
\newcommand{\Db}{\bm {D}}

\RequirePackage{fix-cm}
\begin{document}


\title{Randomized Functional Sparse Tucker Tensor for Compression and Fast Visualization of Scientific Data}

\author{P. Rai        \and
         H. Kolla \and
         L. Cannada \and
         A. Gorodetsky
}


\institute{P.Rai, H. Kolla, L. Cannada \at
		Sandia National Laboratories, Livermore, USA\\
              Correspondence \email{pmrai@sandia.gov}           
           \and
           A. Gorodetsky \at
           Aerospace Engineering, University of Michigan, Ann Arbor, USA 
}

\date{Received: date / Accepted: date}

%


\maketitle

\begin{abstract}
We propose a strategy to compress and store large volumes of scientific data represented on unstructured grids. Approaches utilizing tensor decompositions for data compression have already been proposed. Here, data on a structured grid is stored as a tensor which is then subjected to appropriate decomposition in suitable tensor formats. Such decompositions are based on generalization of singular value decomposition to tensors and capture essential features in the data with storage cost lower by orders of magnitude. However, tensor based data compression is limited by the fact that one can only consider scientific data represented on structured grids. In case of data on unstructured meshes, we propose to consider data as realizations of a function that is based on functional view of the tensor thus avoiding such limitations. The key is to efficiently estimate the parameters of the function whose complexity is small compared to the cardinality of the dataset (otherwise there is no compression). Here, we introduce the set of functional sparse Tucker tensors and propose a method to construct approximation in this set such that the resulting compact functional tensor can be rapidly evaluated to recover the original data. The compression procedure consists of three steps. In the first step, we consider a fraction of the original dataset for interpolation on a structured grid followed by sequentially truncated higher order singular value decomposition to get a compressed version of the interpolated data. We then fit singular vectors on a set of functional basis using sparse approximation to obtain corresponding functional sparse Tucker tensor representation. Finally, we re-evaluate the coefficients of this functional tensor using randomized least squares at a reduced computational complexity. This strategy leads to compression ratio of orders of magnitude on combustion simulation datasets. 
\end{abstract}


\section{Introduction}
Functional tensors are based on interpretation of high dimensional functions as tensors and their decomposition in several tensor formats as particular approximations. Consequently, functional tensors have been studied and applied for sampling based approximation of high dimensional functions in cases where the number of available function evaluations is small. Several functional tensor formats have been studied for various applications e.g. \cite{Doostan:2013} \cite{Chevreuil:2015}\cite{Gorodetsky:2018} \cite{Gorodetsky:2019}. These approaches rely on linearity between the parameters of the low-rank format and the output of the function. Utilizing this multilinear parameterization, they convert the low-rank function approximation to one of low-rank tensor decomposition for the coefficients of a tensor-product basis. 

The novelty of the present paper, in contrast, is aimed at detecting low rank structure in the large volumes of data in order to obtain a low complexity functional tensor representation for a small loss in accuracy. As opposed to high dimensional function approximation using tensors in earlier works, high dimensionality does not come from the number of function parameters but from the number of data points required to be processed in order to obtain a functional tensor form. This functional tensor, stored as a surrogate at a fraction of cost of the original dataset, can be rapidly evaluated to recover accurate approximations of the data. The compressed functional form can act as a preview of the full dataset, which may reside on long-term storage and need not replace the original dataset. In this paper, we consider a particular type of tensor format i.e. functional sparse Tucker representation of the data noting that the ideas presented here can be readily extended to other tensor formats. 

Several compression methods largely focus on compressing local structure with very little loss in precision. Examples of such methods include multivariate volume block data reduction by taking advantage of local multiway structure \cite{Fout:2005}, compression of data in local blocks \cite{lindstrom:2014,Di:2016}. Tensor based methods, in contrast, aim at detecting global structure in the data. It does not process the data in blocks but rather considers the data in its entirety. In this work, in order to take advantage of tensor based compression, we first interpolate the unstructured data on a structured grid followed by its Tucker decomposition \cite{Tucker:1966}. Singular vectors with truncated rank for each mode thus obtained are represented as functions on a suitable basis using least squares with sparsity constraints thus resulting in a functional sparse Tucker representation of the dataset. Finally, to compensate for the effect of interpolation of data for Tucker compression, we re-estimate the components of the core tensor, at a much lower computation complexity than classical approaches, by solving a randomized least squares problem using data in the original dataset.     

The manuscript is organized as follows. We introduce and formalize the notion of functional sparse Tucker tensors in section \ref{func_tucker}. In order to construct approximations in this set, we review least squares with sparse regularization in section \ref{sparse_reg}. We then present our construction algorithm in section \ref{func_tuckermpi} and apply a randomized method for re-estimation of core tensor in section \ref{randls}. Finally, we illustrate our results on combustion simulation datasets in section \ref{appli} with a short conclusion in section \ref{conclusion}.

\section{Functional Sparse Tucker Tensors}
\label{func_tucker}
The key idea in this work is to represent the dataset as realizations of a multivariate function
\begin{align*}
u(y_1,\ldots,y_d) = \sum_{i_1=1}^{n_1}\cdots\sum_{i_d=1}^{n_d}\beta_{i_1,\ldots,i_d}\phi_{i_1}^{(1)}(y_1)\cdots \phi_{i_d}^{(d)}(y_d),
\end{align*} 
where $\phi^{(k)}_{i_k}, 1\leq k\leq d$ are basis functions (e.g. polynomials, wavelets...). The number of expansion coefficients $\beta_{i_1,\ldots,i_d}$ are $\prod_{k=1}^{d}n_k$ thus manifesting the curse of dimensionality if $n_k$ or $d$ or both are large. In such cases, we instead represent the data as realizations of a Tucker low rank approximation $\tilde{u}$ of $u$ where
\begin{align}
\tilde{u}(y_1,\ldots,y_d) = \sum_{j_1=1}^{r_1}\cdots\sum_{j_d=1}^{r_d}\alpha_{j_1,\ldots, j_d}w_{j_1}^{(1)}(y_1)\cdots w_{j_d}^{(d)}(y_d).
\label{eq:func_tuck}
\end{align}
Storage of $\tilde{u}$ in \eqref{eq:func_tuck} require $\prod_{k=1}^{d}r_k$ coefficients and $\sum_{k=1}^dn_kr_k$ expansion coefficients of $w_{j_k}^{(k)}(y_k),1\leq j_k\leq r_k, 1\leq k\leq d$ such that
\begin{align}
w_{j_k}^{(k)}(y_k) = \sum_{i_k=1}^{n_k} w^k_{i_k,j_k} \phi^{(k)}_{i_k}(y_k).
\label{eq:sparse}
\end{align}
In additional, in order to gain advantage from sparsity based regularization, we also constraint the number of non zero coefficients in \eqref{eq:sparse}.
In the following, we formalize the notion of functional sparse Tucker tensors.
\par
We introduce approximation spaces $\Sc^k_{n_k}$ with orthonormal basis $\{\phi_{j}^{(k)}\}_{j=1}^{n_k}$, such that
\begin{align*}
\Sc_{n_k}^k &= 
\left\{v^{(k)}(y_k) = \sum_{j=1}^{n_k} v_j^{k}
\phi_{j}^{(k)}(y_k);v_j^k\in\Rbb\right\} \\
&= \left\{v^{(k)}(y_k) 
= \phib^{(k)}(y_{k})^{T} \v^{(k)}
;\v^{(k)}\in\Rbb^{n_{k}}\right\},
\end{align*}
where $\v^{(k)}$ denotes the vector of coefficients of $v^{(k)}$ and 
where $\phib^{(k)}=(\phi_1^{(k)},\ldots ,\phi_{n_k}^{(k)})^T$ denotes the vector of basis functions.
 An approximation space $\Sc_{\boldsymbol{n}}$ is then obtained by tensorization of approximation spaces $\Sc^k_{n_k} $:
\begin{align*}
\Sc_{\boldsymbol{n}} &= \Sc_{n_1}^1\otimes \hdots \otimes \Sc_{n_d}^d
=\left\{v = \sum_{i \in I_{\boldsymbol n}} 
v_{i} \phi_{i} \; ; \; v_{ i} \in\Rbb\right\},
\end{align*}
where $I_{\boldsymbol n} =  \times_{k=1}^d
\{1 \hdots n_k\}$ and 
$\phi_{ i}(y)=(\phi_{i_1}^{(1)} \otimes \hdots \otimes
\phi_{i_r}^{(d)})(y_1,\hdots,y_d) = \phi_{i_1}^{(1)}(y_1)
\hdots  \phi_{i_r}^{(d)}(y_r)$. An element
$v = \sum_{i} v_{i} \phi_{i} \in\Sc_{\boldsymbol{n}}$ can be identified with the algebraic
tensor $\mathbf{v}\in  \Rbb^{n_1}\otimes
\hdots\otimes \Rbb^{n_d}$ such that $ (\v)_{ i} = v_{ i}$. 
Denoting $\phib(y)=\phib^{(1)}(y_1)\otimes\ldots\otimes\phib^{(d)}(y_r) \in \Rbb^{n_1}\otimes
\hdots\otimes \Rbb^{n_d} $, we have the identification $\Sc_{\boldsymbol{n}} \simeq \Rbb^{n_1}\otimes
\hdots\otimes \Rbb^{n_d}$ with 
$$
\Sc_{\boldsymbol{n}} = \left\{ v(y) = \langle \boldsymbol{\phi}(y),\v \rangle ; \v\in \Rbb^{n_1}\otimes
\hdots\otimes \Rbb^{n_d}\right\},
$$
where $\langle\cdot,\cdot\rangle$ denotes the canonical inner product in $\Rbb^{n_1}\otimes
\hdots\otimes \Rbb^{n_d}$.

Here, we suppose that the approximation space $\Sc_{\boldsymbol n}$ is sufficiently rich to allow accurate representations of a
large class of functions (e.g. by choosing polynomial spaces with
high degree, wavelets with high resolution...). We now introduce the set of functional sparse tensors.

Let $\Rc_1 $ denote the set 
of (elementary) rank-one tensors in $\Sc_{\boldsymbol{n}} =
\Sc^1_{n_1}\otimes \hdots \otimes \Sc^d_{n_d}$, defined by 
\begin{align*}
\Rc_1 &= \left\{w(y) = \left(\otimes_{k=1}^d w^{(k)}\right)(y) =
\prod_{k=1}^d w^{(k)}(y_k) \; ; \; w^{(k)} \in \Sc_{n_k}^k\right\},
\end{align*}
or equivalently by
\begin{align*}
\Rc_1 = \left\{w(y) = \langle\phib(y),{\w^{(1)}}\otimes\ldots \otimes{\w^{(d)}}\rangle;\w^{(k)}\in\Real^{n_k}\right\},
\end{align*}
where $\phib(y) = \phib^{(1)}(y_1)\otimes\hdots\otimes \phib^{(d)}(y_d)$, with $\phib^{(k)}=(\phi_1^{(k)},\ldots ,\phi_{n_k}^{(k)})^T$ the vector of basis functions of  $\Sc^{k}_{n_k}$, and where 
$\w^{(k)}=(w_1^k,\hdots,w_{n_k}^k)^T$ is the set of coefficients of $w^{(k)} $ in the basis of $ \Sc_{n_k}^k$, that means 
$w^{(k)}(y_k) = \sum_{i=1}^{n_k} w_i^k\phi_i^{(k)}(y_k) $. 
Correspondingly, we define  $\boldsymbol{m}$-sparse rank-one subset defined as 
\begin{align*}
\Rc_1^{\boldsymbol{m}\text{-sparse}} &= \Big\{w(y) = \langle \phib(y),{\w^{(1)}}\otimes\ldots \otimes{\w^{(d)}}\rangle;
&\w^{(k)}\in\Real^{n_k},\Vert \w^{(k)} \Vert_0 \leq m_k\Big\}
\end{align*}
with effective dimension $\sum_{k=1}^d m_k\ll \sum_{k=1}^d n_k$ (here we only count the values of the non-zero coefficients and not the integers indicating their locations). However performing 
least-squares approximation in this set may not be computationally tractable. We thus introduce a convex relaxation of the $\ell_0$-``norm''
to define the subset $\Rc_1^\gammab$ of $\Rc_1$ defined as
\begin{align*}
\Rc_1^{{\gammab}} &=\Big\{w(y) = \langle \phib(y),{\w^{(1)}}\otimes\ldots \otimes{\w^{(d)}}\rangle;
&\w^{(k)}\in\Real^{n_k},\Vert \w^{(k)} \Vert_1 \leq \gamma_k\Big\},
\end{align*}
where the set of parameters $(\w^{(1)},\hdots,\w^{(d)})$ is now searched in a convex subset of $\Rbb^{n_1}\times \hdots \times \Rbb^{n_d}$. 
\\
Finally, we introduce the set of functional Tucker tensors with multilinear Tucker rank $\r = (r_1,\ldots,r_d)$  
\begin{align*}
\Tc_{\r} = \left\{ v = \sum_{j_1=1}^{r_1}\cdots\sum_{j_d = 1}^{r_d}\alpha_{j_1,\ldots,j_d} w_{j_i,\ldots,j_d} \;; w_{j_1,\ldots,j_d}\in\Rc_1  \right\} 
\end{align*} 
and the corresponding sparse subset
\begin{align*}
\Tc_{\r}^{\boldsymbol{\gamma}}  = \left\{ v = \sum_{j_1=1}^{r_1}\cdots\sum_{j_d = 1}^{r_d}\alpha_{j_1,\ldots,j_d} w_{j_i,\ldots,j_d} \;; w_{j_1,\ldots,j_d}\in\Rc_1^{\gamma}  \right\}.
\end{align*}
In the following, we propose algorithms for the
construction of approximations in tensor subsets $\Tc_{\r}^{\boldsymbol{\gamma}} $ which requires sparse approximation of functions $w^{k}_{j_k}(y_k)$. For this purpose, we use least squares with sparse regularization as described in the next section.

\section{Least squares with sparse regularization}\label{sparse_reg}
 A sparse function is one that can be represented using few non zero terms when expanded on a suitable
basis. In general, a successful reconstruction of sparse
solution vector depends on  sufficient sparsity of the
coefficient vector and additional properties (incoherence) depending on the samples and of the chosen basis (see \cite{CAN06b,DON06}).
More precisely, an approximation $ \sum_{i=1}^P v_i\phi_i(y)$ of a
function $u(y)$ is considered as sparse on a particular basis
$\{\phi_i(y)\}_{i=1}^P$ if it admits a good approximation with only a few 
non zero coefficients. Under certain conditions,
a sparse approximation can be computed accurately using only
$Q \ll P$ samples of $u(y)$ via sparse regularization.
Given the random samples $\mathbf{z}\in\Rbb^Q$ of the function $u(y)$ at sample points $\{y^q\}_{q=1}^Q$, a best $m$-sparse (or $m$-term) approximation 
of $u$ can be ideally obtained by
solving the constrained optimization problem 
\begin{align}
\min_{\v\in\Rbb^P}  \|{ \mathbf{z}-\Phib\v }\|_2^2 \quad \mathrm{subject\;to}
\quad \Vert{\v}\Vert_0 \le m , \label{eq:min_P0delta}
\end{align}
 where $\norm{\v}_0 = \#\{i\in
\{1,\ldots,P\}\;:\;v_i\neq 0\}$ is the so called $\ell_{0}$-``norm'' of $\v$ which gives the number of non zero components
of $\v$ and and $\Phib \in \Real^{Q\times P}$
the matrix with components $ (\Phib)_{q,i} = \phi_i(y^q)$.
Problem \eqref{eq:min_P0delta} is a combinatorial optimization problem which is NP hard to solve. Under certain assumptions, problem
\eqref{eq:min_P0delta} can be reasonably well approximated by
the  following constrained optimization problem which introduces a convex relaxation of the $\ell_{0}$-``norm'':
\begin{align}
\min_{\v\in\Rbb^P} \|{ \mathbf{z}-\Phib\v }\|_2^2 \quad \mathrm{subject\;to} \quad  \|{\v}\|_1
 \le \delta, 
\label{eq:min_P1delta}
\end{align}
where $\norm{\v}_1 = \sum_{i=1}^P |v_i|$ is the $\ell_1$-norm of $\v$. Since the $\ell_2$ and $\ell_1$-norms are convex, we can equivalently consider the following convex optimization problem, known as Lasso \cite{TIB96} or  basis pursuit \cite{CHE99}:
\begin{align}
\min_{\v \in \Rbb^P} \|{\mathbf{z}-\Phib\v}\|_2^2+\lambda\|\v\|_1,
\label{eq:min_P1lambda}
\end{align}
where $\lambda>0$ corresponds to  a Lagrange multiplier whose value is related to 
$\delta$. Problem \eqref{eq:min_P1lambda} appears as a regularized least-squares problem. The $\ell_{1}$-norm is a sparsity inducing regularization function in the sense that the solution $\v$ of  \eqref{eq:min_P1lambda} may contain components which are exactly zero.  Several optimization algorithms have been
proposed for solving \eqref{eq:min_P1lambda}
(see \cite{BAC12}).  In this paper, we use the Lasso modified least angle 
regression algorithm (see LARS presented in \cite{Efron2004}) and fast leave-one-out cross validation error estimate \cite{CAW04} for optimal sparse solution (corresponding to regularization parameter $\lambda$) which relies on the use of the Sherman-Morrison-Woodbury formula (see  \cite{BLA11} for its implementation within Lasso modified LARS algorithm). In this work, we have used Lasso modified LARS implementation of SPAMS software \cite{MAI10} for $\ell_1$-regularization.

\section{Functional Sparse Tucker Using TuckerMPI}
\label{func_tuckermpi}
\subsection{Interpolation on structured grid}
Representation of the dataset in functional sparse Tucker format defined in section \ref{func_tucker} requires estimation of the core tensor $\alphab$ and univariate functions $w_j^{(k)}(y_k)$. If the dataset is available on a structured grid, it can be stored as a tensor $\Uc$ which can then be decomposed in Tucker format 
\begin{align*}
\Uc \approx \tilde{\Uc} = \alphab\times_1 W^{(1)}\times_2 W^{(2)}\cdots\times_{d} W^{(d)},
\end{align*}
where $\times_k$ is mode $k$ product of $\Uc$ with a factor matrix $W^{(k)}\in \Rbb^{n_k\times r_k}$. Here, the compression precision is given by 
\begin{align*}
\epsilon = \frac{\Vert \Uc - \tilde{\Uc}\Vert_F}{\Vert\Uc\Vert_F},
\end{align*}
where $\Vert\cdot\Vert_F$ is the Frobenius norm. Since the dataset considered is unstructured, we propose to interpolate the data on a structured grid. Let us denote the grid size in mode $k, 1\leq k\leq d$ as $I_k$ and, for the sake of simplicity, consider that the gird points are equispaced. A structured grid of size $I_1\times I_2\cdots\times I_d$ can thus be obtained. Now, we consider only a small subset of the original dataset, say 10\%, for linear interpolation on this grid and the interpolated data is stored as a tensor which is then decomposed in Tucker format. We use TuckerMPI \cite{TUCKERMPI1}, a parallel C++/MPI software package for compressing distributed data, for this purpose. Note that TuckerMPI is a parallel implementation of the sequentially-truncated HOSVD (ST-HOSVD) \cite{STHOSVD}. We thus obtain factor matrices $W^{(k)}, 1\leq k\leq d$, the columns of which are realizations of univariate functions $w_{j_k}^{(k)}(y_k), 1\leq j\leq r_k$.  
\subsection{Sparse approximation of singular vectors}
We now wish to obtain a functional representation $w_{j_k}^{(k)}(y_k), 1\leq j_k\leq r_k$, of the singular vectors $W_{:,j_k}^{(k)}$ such that $W_{i_k,j_k}^{(k)}$ are evaluations of $w_{j_k}^{(k)}(y_k^{i_k})$ at grid locations $\{y_k^{i_k}\}_{i_k=1}^{I_k}$ along mode $k$. For this purpose, we use least squares with sparse regularization in section \ref{sparse_reg} to obtain coefficients on suitable basis functions. It is well known that singular vectors are more oscillatory (see for e.g. Figure \ref{last_mode}(a)) for higher rank as they capture high frequency phenomenon in the dataset. Thus, choice of basis functions for representation of $w_{j_k}^{(k)}(y_k)$ corresponding to small $j_k$ may not be appropriate for the ones with higher $j_k$. Therefore, in this work, we propose to construct approximation in two spaces $\Pbb_{p}$, where $\Pbb_{p}$ is the space of Legendre polynomials of degree $p$ and $\Wbb_{s,p}$, where $\Wbb_{s,p}$ is the space of multi-resolution wavelets with resolution $s$ and degree $p$. We can then choose the approximation that gives smaller approximation error. We present the overall compression scheme in Algorithm \ref{func_tucker_comp} below.

\begin{algorithm}
\caption{Compression of unstructured data in functional sparse Tucker format}
\label{func_tucker_comp}
\begin{algorithmic}[1]
\REQUIRE Original dataset, interpolation grid $I_k,1\leq k\leq d$, compression precision $\epsilon$
\ENSURE Function sparse Tucker tensor core $\alphab$ and  coefficients of $w^{(k)}_{j_k} (y_k), 1\leq k\leq d, 1\leq j_k\leq r_k$.
\STATE Interpolate the data on structured grid of size $I_1\times \cdots \times I_d$
\STATE Use TuckerMPI to get core tensor $\alphab$ and factor matrices $W^{(k)}$ for given compression precision 
\FOR {$k=1,\ldots,d$}
\FOR {$j_k=1,\ldots,r_k$}
\STATE Approximate $w^{(k)}_{j_k}$ using components of $W^{(k)}_{:,j_k}$ in $\Pbb_p$ and $\Wbb_{s,p}$ and estimate error (See section \ref{sparse_reg})
\STATE Store coefficients of $w^{(k)}_{j_k}$ corresponding to smaller approximation error
\ENDFOR
\ENDFOR
\end{algorithmic}
\end{algorithm}

\section{Re-estimation of core tensor using randomized least squares}
\label{randls}
Algorithm \ref{func_tucker_comp} gives a functional representation of the data in the sparse Tucker format. The data however was interpolated on a structured gird to be able to use TuckerMPI for obtaining the singular vectors. This approach suffers from the limitation that the quality of approximation will depend on the type of interpolation (e.g. linear, non linear interpolation, cardinality of points considered for interpolation etc.) on a structured grid. To overcome this limitation, it is imperative to use the original dataset to re-evaluate some, if not all, parameters of functional sparse Tucker tensor. We implement this idea by re-evaluating the elements of the core tensor using linear least squares. Let us rewrite \eqref{eq:func_tuck} as 
\begin{align*}
u(y)=\tilde{u}(y) = \sum_{\jb=1}^{R} \alpha_{\jb} w_{\jb}(y),
\end{align*}
where $\jb = (j_1,\ldots, j_d)$ such that $\alpha_{\jb} = \alpha_{j_1,\ldots, j_d}$ and $w_{\jb}(y) = \prod_{k=1}^{d}w_{j_k}^{(k)}(y_k)$ and $R=r_1\times\cdots \times r_d$. We wish to solve the regression problem
\begin{align}
\hat{\alphab} = \min_{\alphab\in\Rbb^{R}} \|{ \mathbf{u}-\Wb\alphab}\|_2^2,
\label{rand_core}
\end{align}
where $(\u)_{q} = u(y^q)$, $\Wb \in \Rbb^{Q\times R}$ is the matrix such that $(\Wb)_{q,{\jb}} = w_{\jb}(y^q)$ and $\alphab \in \Rbb^{R}$ are components of the core tensor reshaped as a vector. Clearly, $\Wb$ is overdetermined $(Q\gg R),$ and hence its computation and storage may be prohibitive in the considered setting. However, sketching and randomized methods have been used successfully to solve big linear least squares problems at a much smaller computation cost \cite{ROKHLIN08,DRINEAS11,AVRON10,KANNAN17} and is ideal for application in this setting. 

The key here is to transform \eqref{rand_core} using random projection $\Mb\in \Rbb^{S\times R}, S\ll Q,$ such that an exact solution to $\mathrm{min}_{\alphab}\|{ \Mb\mathbf{u}-\Mb\Wb\alphab}\|_2^2$ is an approximate solution to the original problem \eqref{rand_core} \cite{WOODRUFF14}. Approaches to solve randomized least squares problems are based on the idea of leverage scores. The leverage score of rows of $\Wb$ is the norm of the rows of its left singular vectors and corresponds, in some sense, to the importance of that row in constructing its column-space. One can then solve the randomized least squares problem by sampling rows of $\Wb$ weighted according to the distribution of the leverage scores. One drawback of this method is that leverage scores have to be estimated from singular value decomposition of $\Wb$ which may be computationally expensive when $R$ is large. Therefore, we follow the approach in \cite{BATTAGLINO18}, where we mix $\Wb$ with the intention of evenly distributing leverage scores across all rows in such a way that one can sample rows uniformly. Note that this mixing strategy relates to a more general class of transformations that rely on approximation quality guarantees provided by the Johnson-Lindenstrauss Lemma \cite{JOHNSON84}. This lemma specifies a class of random projections that preserve the distances between all pairs of vectors with reasonable accuracy in the projected subspace.

We briefly summarize steps for our problem. Firstly, we premultiply $\Wb$ by a diagonal matrix $\Db\in \Rbb^{Q\times Q}$ with random +1/-1 to spread out the signal in frequency domain \cite{AILON06}. This is equivalent to flipping the sign of each row of $\Wb$ with probability 1/2. We then apply a fast mixing operation, here a fast Fourier transformation, which has the effect of mixing information across every element of a vector. At the end of this step, the leverage scores of the resulting matrix are concentrated about a small value. Finally, we sample S rows of this matrix with uniform probability. These steps define the random projection $\Mb$. Algorithm \ref{alg:randlsq} outlines the steps in re-estimation of the core tensor.

\begin{algorithm}
\caption{Re-evaluation of the core tensor using Randomized Least Squares}
\label{alg:randlsq}
\begin{algorithmic}[1]
\REQUIRE Original dataset $\u$, Core tensor $\alpha$, Functional singular vectors $w_{j_k}^{(k)}, 1\leq k\leq d, 1\leq j_k\leq r_k$
\ENSURE Re-evaluated Function sparse Tucker tensor core $\hat{\alphab}$
\STATE Construct $\Wb$ by evaluating $w_{\jb}$ at samples in the dataset
\STATE Multipy $\Wb$ and $\u$ with diagonal matrix $\Db$ 
\STATE Apply a fast Fourier Transformation i.e. $\Fc\Db\Wb$ and $\Fc\Db\u$
\STATE Sample $S>R$ rows uniformly, i.e. $\Sb\Fc\Db\Wb$ and $\Sb\Fc\Db\u$, where $\Sb$ is a sampling matrix
\STATE Solve $\hat{\alphab} = \mathrm{min}_{\alphab} \Vert \Sb\Fc\Db\Wb - \Sb\Fc\Db\u\Vert_2^2$
\end{algorithmic}
\end{algorithm}

\section{Illustration}
\label{appli}
We apply our compression strategy on a data set pertaining to a direct numerical simulation (DNS) of turbulent combustion. A ``statistically planar'' (SP) premixed flame \cite{SP_flame} stabilized in homogeneous isotropic turbulence is simulated using the massively parallel DNS code S3D \cite{S3D}. A premixed mixture of methane and air establish a flame that remains statistically planar and stationary in an oncoming turbulent flow. The combustion chemistry is described  using a chemical mechanism containing six chemical species. Accordingly, at each point in the spatial grid and time the solution vector contains eleven dependent variables describing the full thermo-chemical state of the flame. The data set is mapped onto a 3-dimensional structured grid comprising 500 grid points in each spatial dimension, and a total of 400 time snapshots are considered. For the illustration of the method, to follow, we consider two
variants of this fundamentally 4-dimensional data set. In the first case, henceforth referred to as SP3D, we consider that the data belongs to a three dimensional space, consisting of two spatial axis and one time axis. The total number of data points in this set is $7.5\times10^7$ with total storage cost of 0.6 gigabytes for double precision. The second case considers a 4 order tensor, SP4D, which also considers the third spatial axis, in addition to the ones in SP3D. The total storage cost of data in this case is 300 gigabytes with  $3.75\times 10^{10}$ data points. In the following, we illustrate results of SP3D case, and mention that a similar illustrations can be obtained for SP4D.

\begin{figure}[h!]
 \centering
\includegraphics[width=0.6\linewidth]{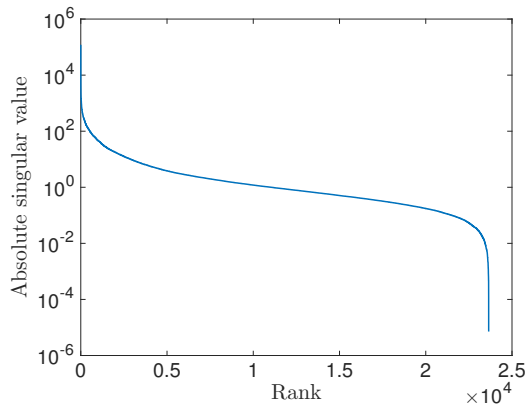}
\caption{Decay of singular values i.e. absolute value of components of core tensor $\alphab$ v/s rank of SP3D with TuckerMPI precision of $1.0\times 10^{-4}$. The rank on horizontal axis is converted to canonical rank by sorting the singular values in descending order.}
\label{sing_vals}
\end{figure}

In case of SP3D, we interpolate the data on a structured grid of size $500\times 500\times 300$ using only 10\% of the data in the original set and decompose the resulting tensor in Tucker format using TuckerMPI. Figure \ref{sing_vals} shows the decay in the absolute value of the components of the core tensor $\alphab$ versus rank (multilinear Tucker rank on horizontal axis is converted to canonical rank) of SP3D for Tucker decomposition precision of $1.0\times10^{-4}$. We clearly see that there is a fast decay in the singular values thus indicating strong scope for compressibility of this dataset. Table \ref{tab:tuckerranks} summarizes the interpolation parameters and multilinear Tucker ranks i.e. size of the core tensor thus obtained for different decomposition precisions for the two datasets.

\begin{table}[h!]
 \centering
  \caption{Specification of interpolation grid size and Tucker ranks obtained for different precision using TuckerMPI for two test cases}
  \label{tab:tuckerranks}
  \begin{tabular}{cccc}
    \hline
    Dataset &Interpolation grid size&Precision $(\epsilon)$&Size of core tensor\\
    \hline
    \multirow{ 2}{*} {SP3D}& \multirow{ 2}{*}{$500\times500\times300$}& $1.0\times 10^{-2}$& $25 \times 24\times 8$\\
    && $1.0\times 10^{-4}$& $57\times 50 \times 17$\\
    {SP4D}& $500\times500\times500\times300$& $1.0\times 10^{-2}$& $30 \times 38 \times 5\times11$\\
  \hline
\end{tabular}
\end{table}

\begin{figure}[h!]
\centering
\subfigure[]{\includegraphics[scale=.30]{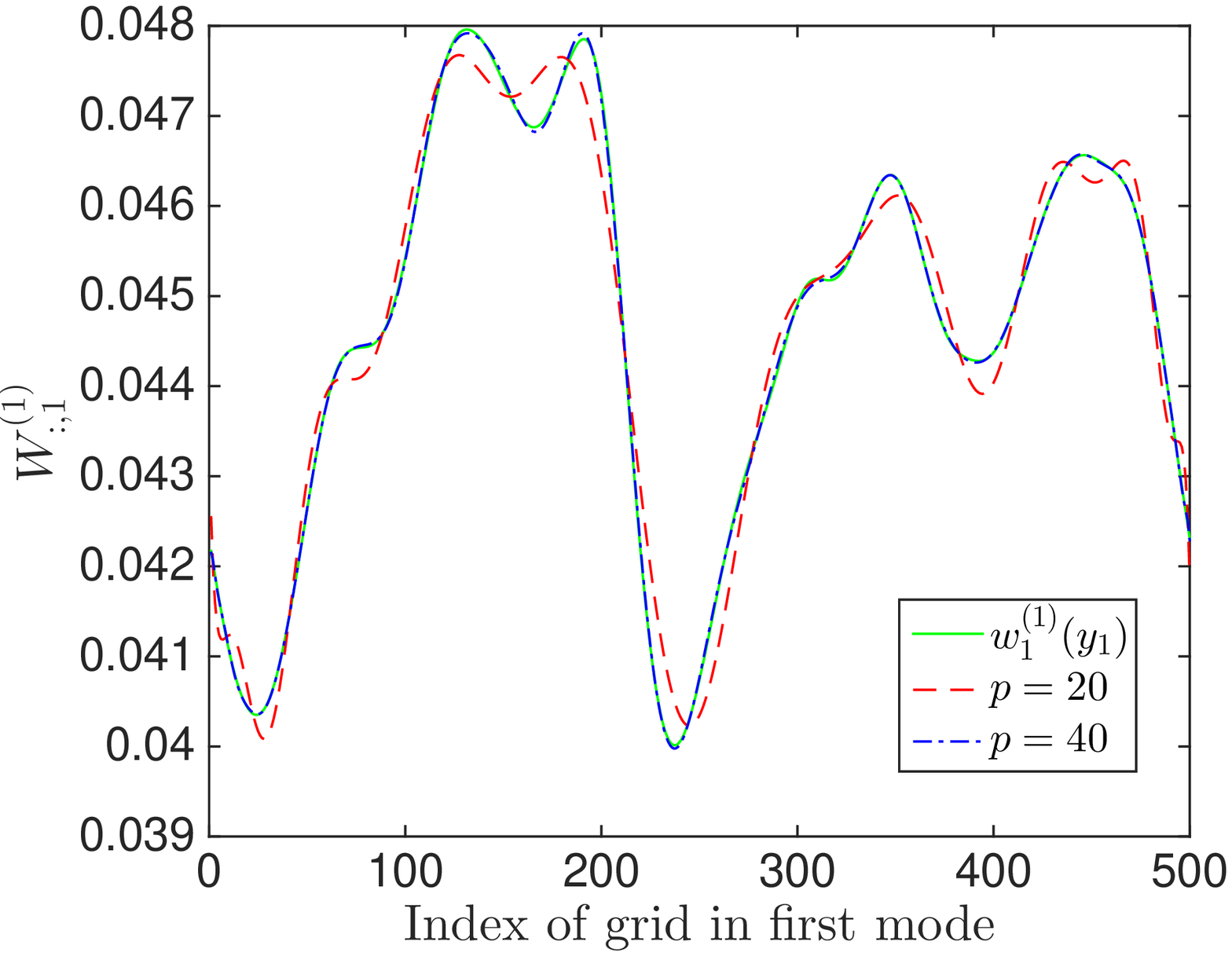}}
\subfigure[]{\includegraphics[scale=.30]{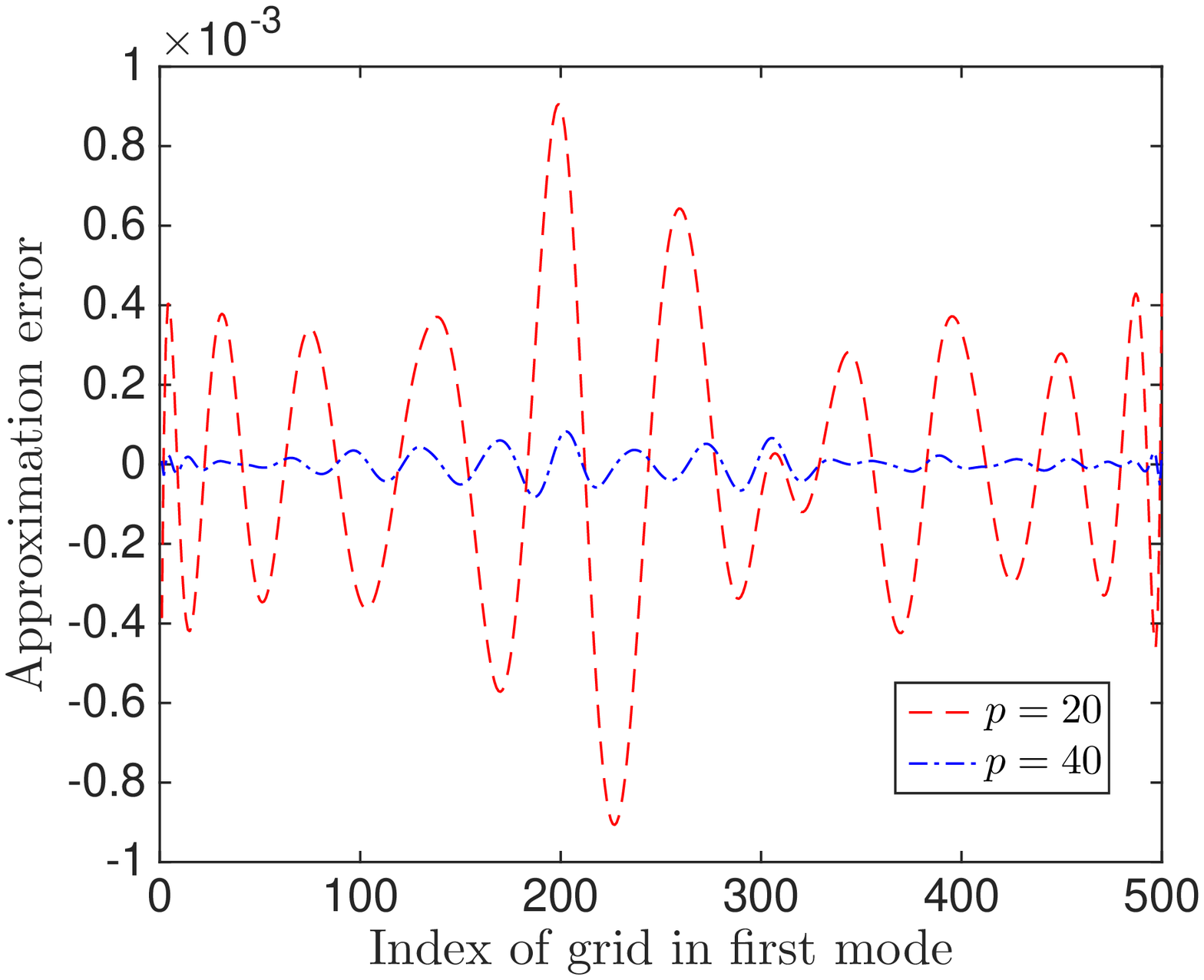}}
\caption{(a) Approximation of $w^{(1)}_1(y_1)$ using least squares with $\ell_1$ regularization from data points as components of $W^{(1)}_{:,1}$ in the approximation space of Legendre polynomials of degree $p = 20$ and $p = 40$. (b) Point wise approximation error v/s grid index of the two approximations in (a)}
\label{first_mode}
\end{figure}

We now consider functional approximations of singular vectors along the first mode. Figure \ref{first_mode}(a) shows first singular vector $W^{(1)}_{:,1}$ and its corresponding functional approximations in $\Pbb_{20}$ and $\Pbb_{40}$. For better illustration, the corresponding approximation errors are plotted in Figure \ref{first_mode}(b). We find that a sufficiently rich approximation space is necessary for accurate representation of singular vectors as point wise error for $p=40$ is much smaller than with $p=20$. 

\begin{figure}[h!]
\centering
\subfigure[]{\includegraphics[scale=.30]{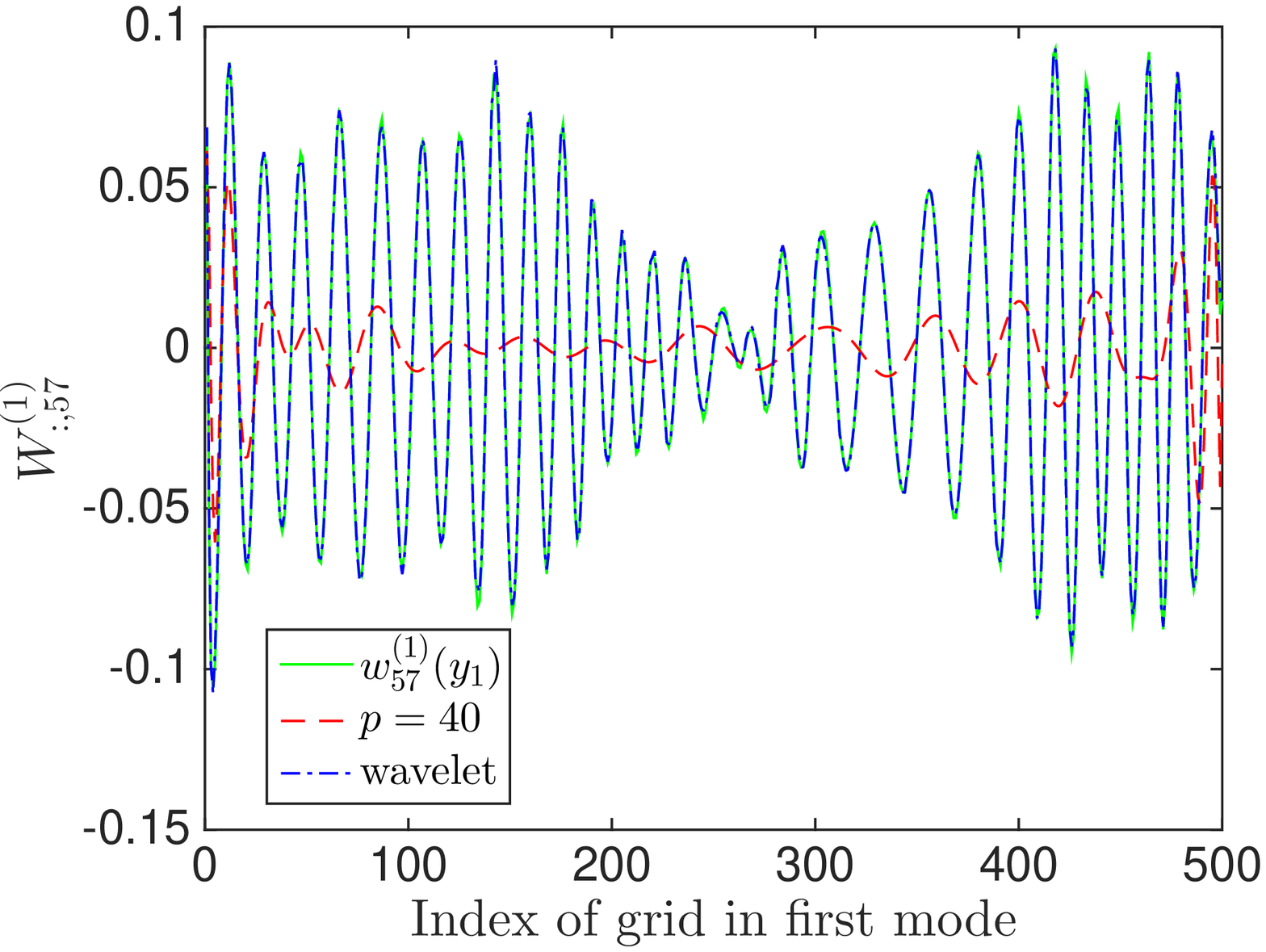}}
\subfigure[]{\includegraphics[scale=.30]{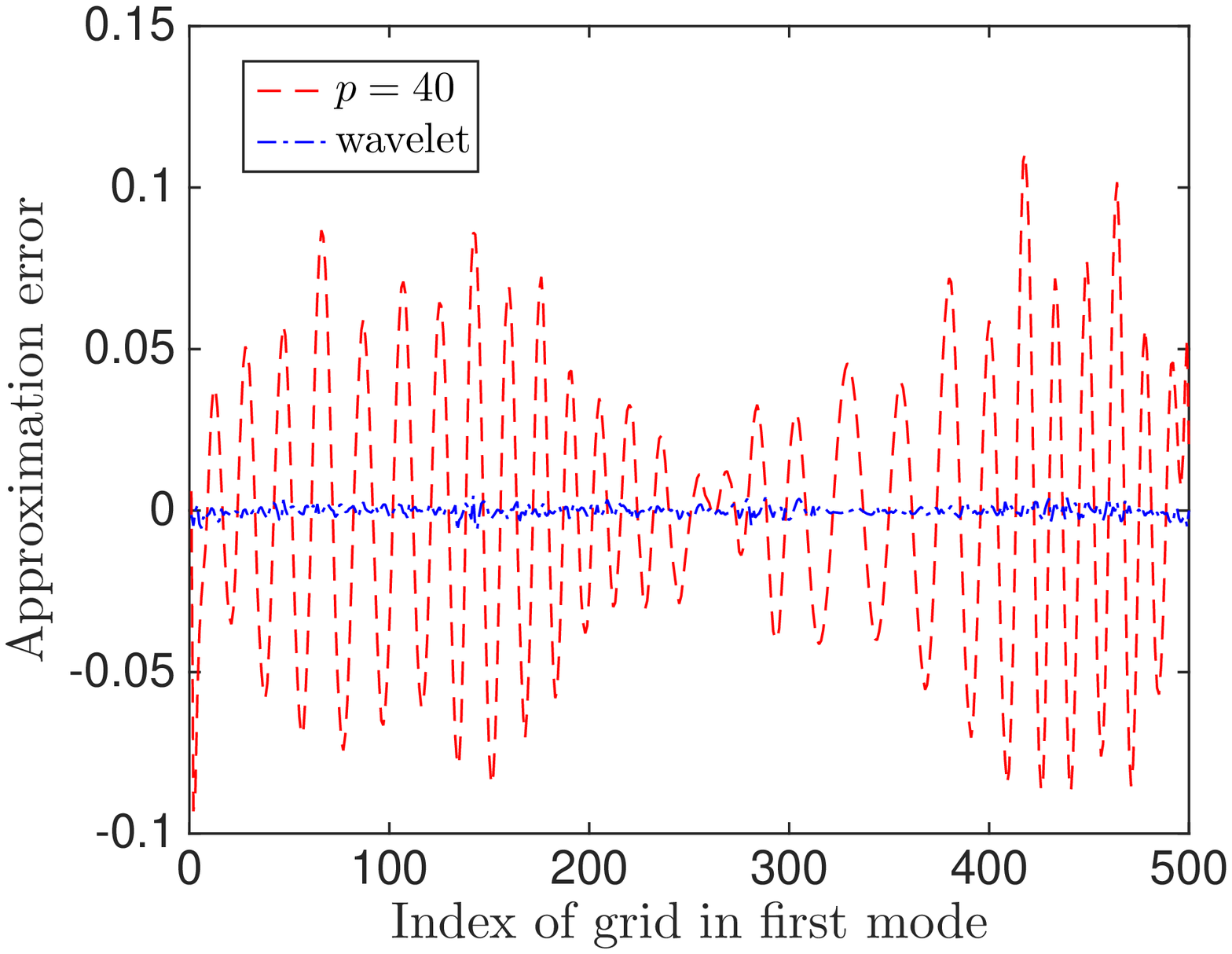}}
\caption{Approximation of $w^{(1)}_{57}(y_1)$ using least squares with $\ell_1$ regularization from data points as components of $W^{(1)}_{:,57}$ in the approximation space of Legendre polynomials of degree $p = 20$  and wavelets with resolution level 5 and degree 3. (b) Point wise approximation error v/s grid index of the two approximations in (a)}
\label{last_mode}
\end{figure}

Figure \ref{last_mode}(a) and (b) show similar plots for the last singular vector in the first mode i.e. $W^{(1)}_{:,57}$ in approximation spaces $\Pbb_{40}$ and $\Wbb_{3,5}$. We clearly see that, in this case, a multi-resolution wavelet basis is essential to get an accurate functional representation, although point wise approximation error is high as compared to the first singular vector. This is because the singular vectors corresponding to higher ranks capture high frequency signals in the dataset and hence require functional basis with higher resolution.

\begin{figure}[h!]
\centering
\subfigure[]{\includegraphics[scale=.28]{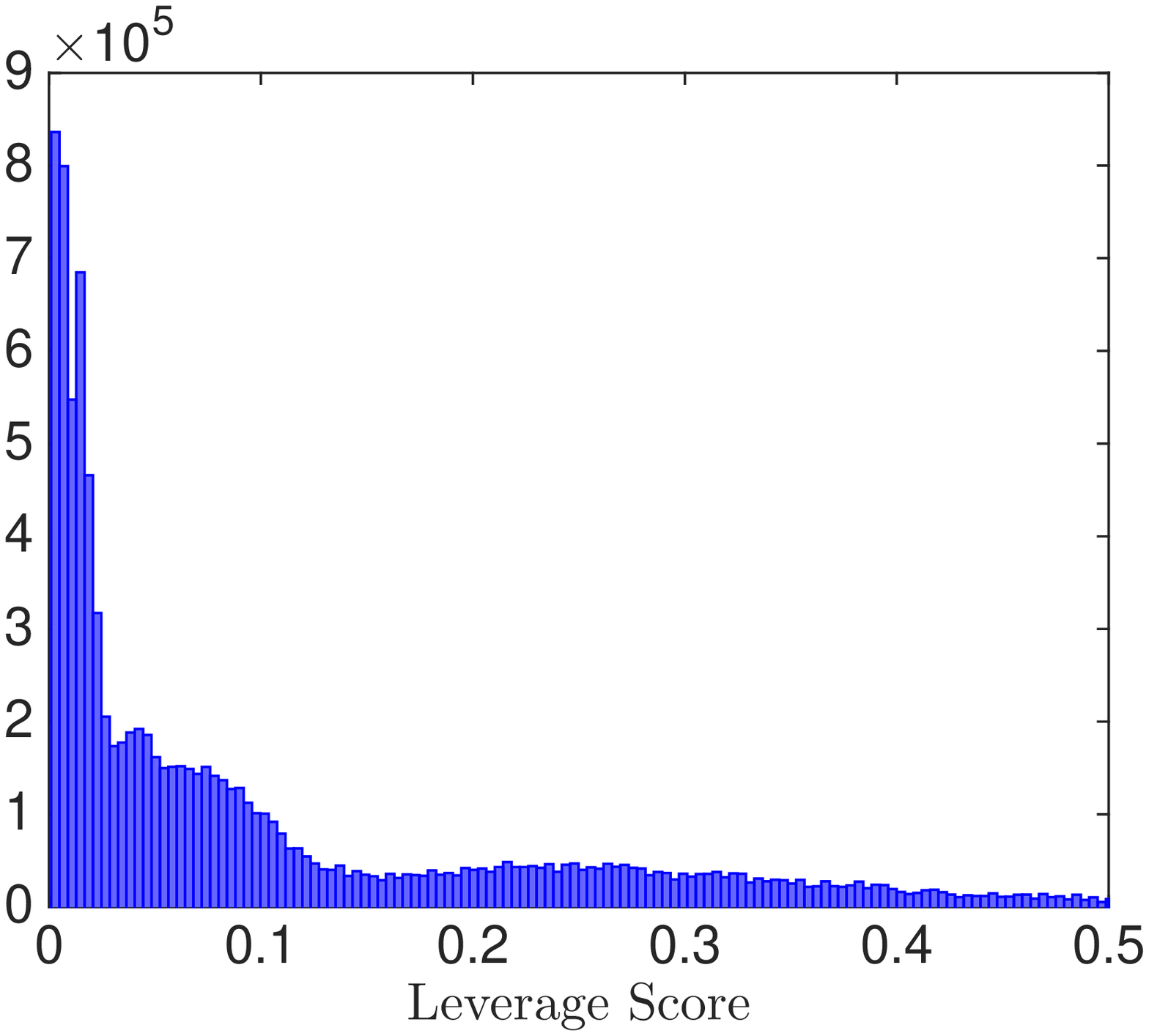}}
\subfigure[]{\includegraphics[scale=.28]{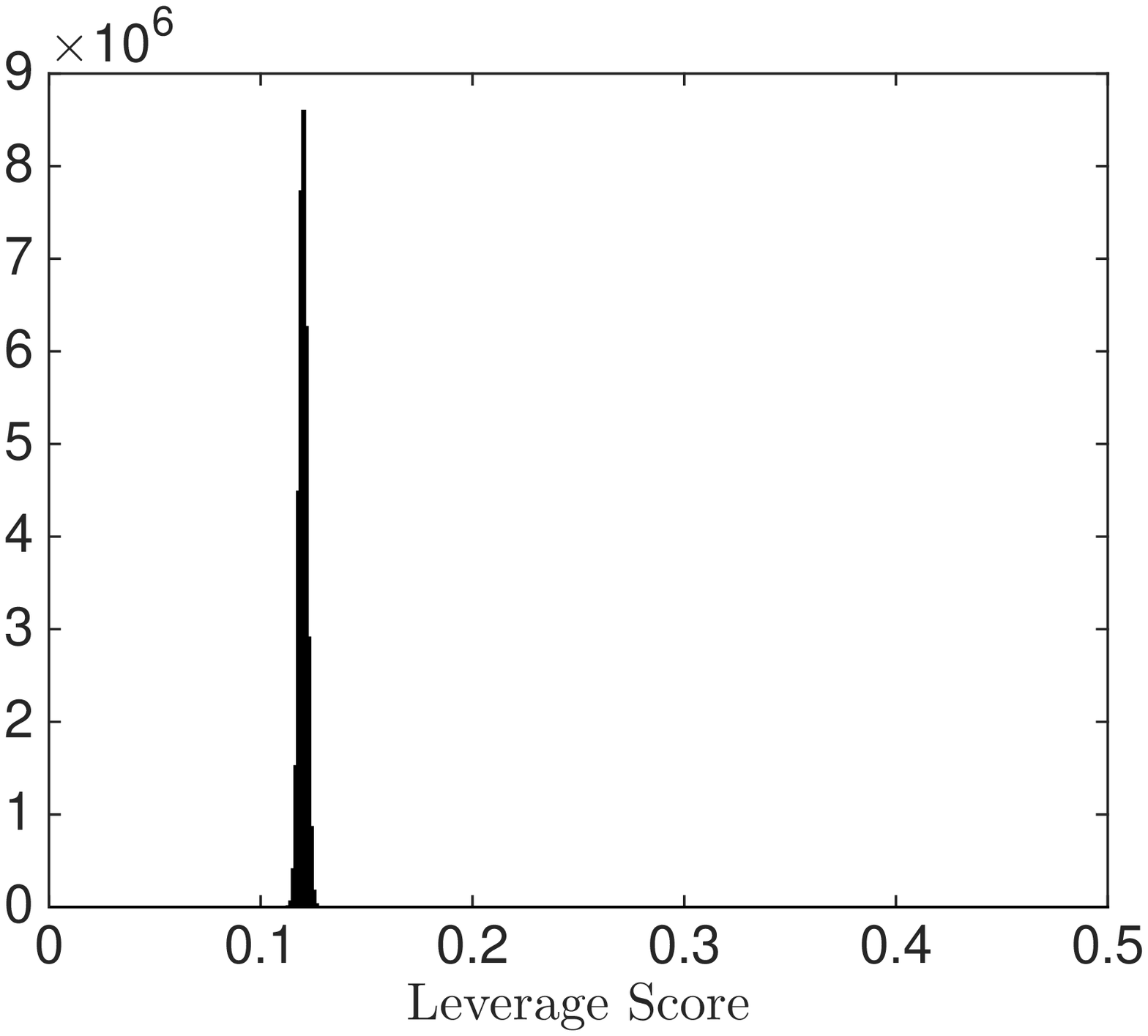}}
\caption{Histogram plot of leverage scores of the measurement matrix $\Wb$ in (a) before mixing and (b) after mixing (Step 3. of Algorithm \ref{alg:randlsq})}
\label{fig:leverage_scores}
\end{figure}

\begin{figure}[h!]
 \centering
\includegraphics[width=0.5\linewidth]{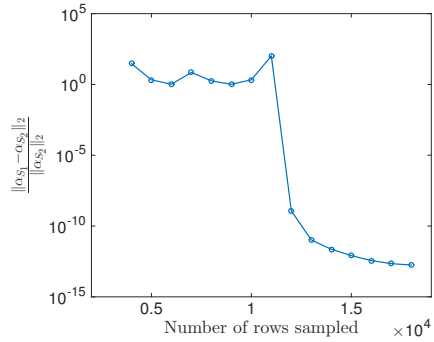}
\caption{Self convergence plot of $\alphab$. Horizontal axis shows the number of rows $S$ sampled in Step 4. of algorithm \ref{alg:randlsq}. Vertical axis shows the change in $\alphab$ (relative norm) estimated by sampling $S_1$ and $S_2$ rows ($S_2 = S_1+10^3, S_1\in\{3000,4000,\ldots,17000\}$).}
\label{randselfconv}
\end{figure}

We now illustrate results related to re-estimation of core tensor using randomized least squares. In Figure \ref{fig:leverage_scores}, we show the distribution of the leverage score of the measurement matrix $\Wb$ (see section \ref{randls}) before and after mixing operation (Step 3. of algorithm \ref{alg:randlsq}). We find that leverage scores, although skewed, have non negligible mass in the range of [0.2,0.5]. On the other hand, after application of mixing operation, leverage scores are concentrated around a small value (0.12). One can thus sample the desired number of rows uniformly to reduce the size of problem \eqref{rand_core}. In Figure \ref{randselfconv}, we illustrate self convergence of the re-estimated core tensor by measuring the relative norm of change in $\alphab$ (by solving step 5. of algorithm \ref{alg:randlsq}) for $S=S_1$ and $S=S_2$, where $S_2>S_1$. We find two distinct regions in the plot separated by a sharp drop in self convergence error at $S \approx2.5R$. Here $R = 25\times 24\times 8$, i.e. size of the core tensor whose coefficients are being re-estimated. Note that an oversampling factor of 2.5 is orders of magnitude smaller as compared to size of the problem when estimating the core tensor with all points in the dataset.

\begin{figure}[h!]
\centering
\subfigure[]{\includegraphics[angle=-90,scale=.20]{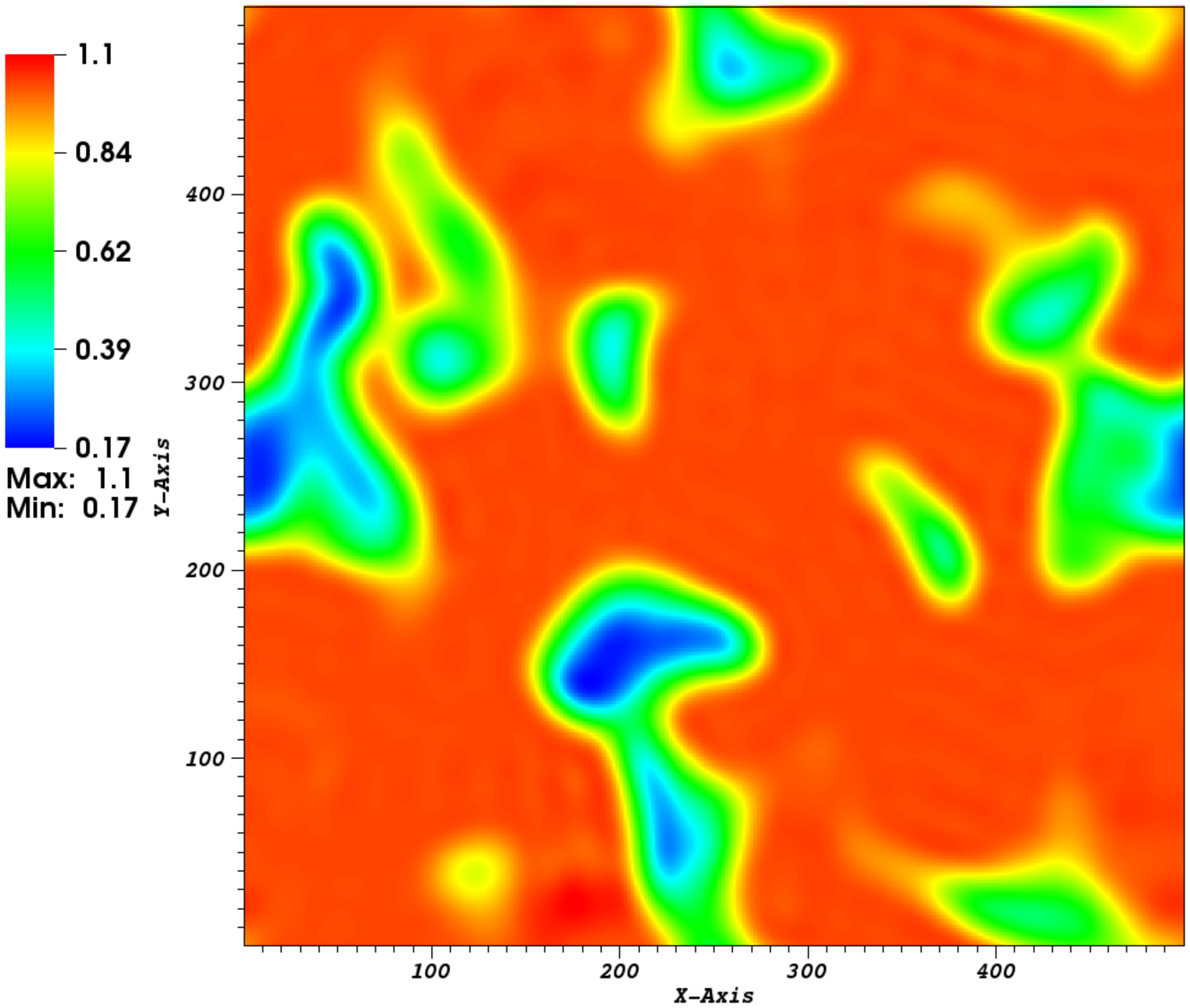}}
\subfigure[]{\includegraphics[angle=-90,scale=.20]{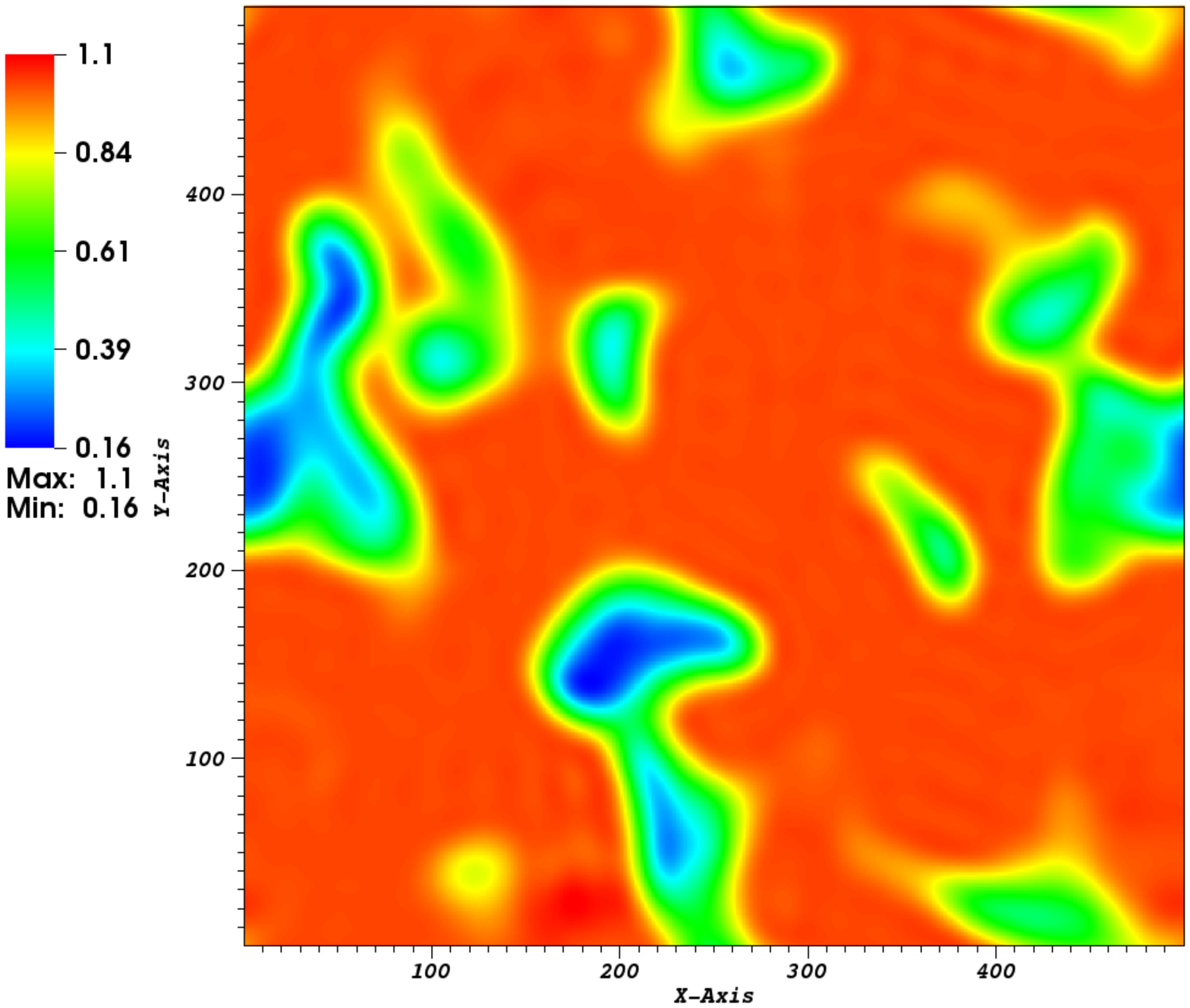}}
\caption{Visualization of a 2D slice of SP4D dataset. (a) Original data and (b) Reconstructed tensor slice obtained from randomized functional sparse Tucker tensor.}
\label{viz}
\end{figure}

Table \ref{tab:compression} shows compression error, compression ratio and storage cost for functional sparse Tucker tensor for both test cases. Finally, Figure \ref{viz} shows visualization of 2D slice of the original dataset obtained from reconstruction of data from sparse functional Tucker tensor. Note that reconstruction of the data from functional tensor is computationally inexpensive as the only computation required is evaluation of the basis functions at points of interest.

\begin{table}[htbp]
 \centering
  \caption{Compression results using functional sparse Tucker tensor}
  \begin{tabular}{ccccc}
    \hline
    Dataset & \makecell{Original \\Storage cost}& Precision&\makecell{Compression \\ratio}&\makecell{Reduced \\Storage cost}\\
    \hline
    \multirow{ 2}{*} {SP3D}& \multirow{ 2}{*}{0.6 GB}&$1.01\times 10^{-2}$& 3879 & 155 KB\\
    &&$1.9\times 10^{-4}$& 936 &640 KB\\
    {SP4D}&300 GB& $1.1\times 10^{-2}$& $4.45\times10^5$& 673 KB \\
  \hline
\end{tabular}
\label{tab:compression}
\end{table}

\section{Conclusion}
\label{conclusion}
We presented a novel technique to compress large volume of data using functional sparse Tucker decomposition. The key idea is to find a sufficiently accurate representation of data in the set of functional Tucker tensors with complexity smaller by orders of magnitude as compared to the size of dataset. In order to achieve this objective, we defined the set of sparse functional Tucker tensors and used existing parallel implementation of Tucker decomposition to construct approximation in this set. The singular vectors are approximated as functions represented on suitable basis using least squares with sparse regularization. The entire compression scheme was tested on datasets obtained from high fidelity combustion modeling simulations. For small loss of accuracy, the proposed strategy results in compression ratio of up to $ 10^3-10^{5}$ for a third order and fourth order dataset respectively.

\section{Acknowledgment}
We would like to acknowledge Tamara Kolda for discussions related to this work. Sandia National Laboratories is a multimission laboratory managed and operated by National Technology and Engineering Solutions of Sandia, LLC., a wholly owned subsidiary of Honeywell International, Inc., for the U.S. Department of Energy's National Nuclear Security Administration.

\end{document}